\documentclass[12pt]{article}
\usepackage{amsfonts,amssymb,latexsym}

\newtheorem{theorem}{Theorem}[section]
\newtheorem{lemma}[theorem]{Lemma}

\newtheorem{proposition}[theorem]{Proposition}

\newenvironment{proof}
{\par\addvspace{0.3cm}\noindent{\rm Proof. }}
{\nopagebreak\mbox{}\hfill $\Box$\par\addvspace{0.25cm}}

\renewcommand{\Im}{\mbox{\rm Im\,}}

\newcommand{\R}{{\mathbb R}}
\newcommand{\C}{{\mathbb C}}
\newcommand{\Z}{{\mathbb Z}}
\newcommand{\T}{{\mathbb T}}
\newcommand{\W}{\mathcal{W}}

\newcommand{\pcea}{PC_{\pm1}^{\mathrm{abs}}}
\newcommand{\Li}{L^\iy(\T)}
\newcommand{\Ta}{\T_{-1}}
\newcommand{\Tb}{\T_1}

\renewcommand{\kappa}{\varkappa}
\renewcommand{\rho}{\varrho}

\newcommand{\be}{\begin{equation}}
\newcommand{\ee}{\end{equation}}

\newcommand{\bqn}{\begin{eqnarray}}
\newcommand{\eqn}{\end{eqnarray}}
\newcommand{\nn}{\nonumber}
\newcommand{\ba}{\begin{array}}
\newcommand{\ea}{\end{array}}

\newcommand{\al}{\alpha}
\newcommand{\iv}{^{-1}}
\newcommand{\iy}{\infty}
\newcommand{\eps}{{\varepsilon}}

\newcommand{\ovl}{\overline}

\newcommand{\diag}{\mathrm{diag}}

\newcommand{\twomat}[1]{\left(\ba{cc}#1\ea\right)}

\newcommand{\ta}{\tilde{a}}
\newcommand{\tb}{\tilde{b}}
\newcommand{\tc}{\tilde{c}}

\begin{document}

\date{}
\title{Dyson's constant in the asymptotics of the Fredholm 
determinant of the sine kernel}
\author{Torsten Ehrhardt
	\thanks{tehrhard@mathematik.tu-chemnitz.de.}\\
	Fakult\"{a}t f\"{u}r Mathematik\\
	Technische Universit\"{a}t Chemnitz\\
	09107 Chemnitz, Germany}
\maketitle

\begin{abstract}
We prove that the asymptotics of the Fredholm determinant of
$I-K_\al$, where $K_\al$ is the integral operator with the sine kernel
$\frac{\sin(x-y)}{\pi(x-y)}$ on the interval $[0,\al]$ is given by
$$
\log\det(I-K_{2\al})=
-\frac{\alpha^2}{2}-\frac{\log\alpha}{4}+\frac{\log 2}{12}+3\zeta'(-1)+o(1),
\qquad\al\to\iy.
$$
This formula was conjectured by Dyson. The first and second order asymptotics
of this formula have already been proved and higher order asymptotics
have also been determined. In this paper we solve the remaining outstanding problem
of identifying the constant (or third order) term.
\end{abstract}


\section{Introduction}

Let $K_{\alpha}$ be the integral operator defined on $L^2[0,\al]$ with the
kernel
\bqn
k(x,y)&=&\frac{\sin(x-y)}{\pi(x-y)}.
\eqn
Dyson \cite{D1} conjectured the following asymptotic formula for the determinant
$\det(I-K_{2\alpha})$,
\bqn\label{f.Dys}
\log \det(I-K_{2\alpha})&=&
-\frac{\alpha^2}{2}-\frac{\log\alpha}{4}+\frac{\log 2}{12}+3\zeta'(-1)+o(1),
\qquad \al\to\iy,
\eqn
and provided heuristic arguments.  Moreover, Jimbo, Miwa, M\^{o}ri and Sato
\cite{JMMS} showed that the function
$$
\sigma(\al)=\al\frac{d}{d\al}\log\det (I-K_\al)
$$
satisfies a Painlev\'e V equation. Widom \cite{W2,W3} determined the highest term 
in the asymptotics of $\sigma(\alpha)$ as $\alpha\to\iy$. Knowing this asymptotics one can 
derive a complete asymptotic expansion for $\sigma(\al)$.
{}From this it follows by integration that the asymptotic expansion of $\det(I-K_{2\al})$
is given by
\bqn
\log \det(I-K_{2\alpha}) &=&
-\frac{\alpha^2}{2}-\frac{\log\alpha}{4}+C+
\sum_{n=1}^N \frac{C_{2n}}{\alpha^{2n}}+O(\alpha^{2N+2}),\quad
\al\to\iy,
\eqn
with effectively computable constants $C_2,C_4,\dots$.
The only remaining problem has been the determination of the constant $C$.
This will be done in the present paper. 
In fact, we will prove the asymptotic formula (\ref{f.Dys}).

Let us remark that asymptotic formulas for the determinants of the sine kernel integral operator
defined on $L^2(\al J)$ where $J$ is a finite union of finite subintervals of $\R$
have been considered, and results were established by Widom \cite{W3} and
by Deift, Its and Zhou \cite{DIZ}.

The determinant $\det(I-K_\alpha)$ appears in random matrix theory \cite{Me}.
It is equal to the probability that in the bulk scaling limit of the 
Gaussian Unitary Ensemble of Hermitian matrices an interval of length $\alpha$ contains no eigenvalues.
For further connections we refer to \cite{BTW} and the literature cited there.


\section{Notation}

Let us first introduce some notation. For a Lebesgue measurable subset $M$
of the real axis $\R$ or of the unit circle $\T=\{t\in\C\;:\;|t|=1\}$, let
$L^p(M)$ ($1\le p<\iy$) stand for the space of all Lebesgue measurable
$p$-integrable complex-valued functions. For $p=\iy$ we denote by
$L^\iy(M)$ the space of all essentially bounded Lebesgue measurable
functions on $M$.

For a function $a\in L^1(\T)$ we introduce the
$n\times n$ Toeplitz and Hankel matrices
\be
T_n(a)=(a_{j-k})_{j,k=0}^{n-1},\qquad
H_n(a)=(a_{j+k+1})_{j,k=0}^{n-1},
\ee
where
$$
a_k =\frac{1}{2\pi}\int_0^{2\pi}a(e^{i\theta})e^{-ik\theta}\, d\theta,\qquad
k\in\Z,
$$
are the Fourier coefficients of $a$.
We also introduce a differently defined $n\times n$ Hankel matrix
\be\label{f.Hn2}
H_n[b]=(b_{j+k+1})_{j,k=0}^{n-1},
\ee
where the numbers $b_k$ are the moments of a function $b\in L^1[-1,1]$,
$$
b_{k}=\frac{1}{\pi}\int_{-1}^1 b(x) (2x)^{k-1}\,dx,\qquad k\ge1.
$$

Given $a\in \Li$ the multiplication operator $M(a)$ acting on $L^2(\T)$ is defined by
\bqn
M(a)&:& f(t)\in L^2(\T)\mapsto a(t)f(t)\in L^2(\T).
\eqn
We denote by $P$ the Riesz projection
$$
P:\sum_{k=-\iy}^\iy f_k t^k\in L^2(\T)\mapsto
\sum_{k=0}^\iy f_k t^k\in L^2(\T)
$$
and by $J$ the flip operator
$$
J:f(t)\in L^2(\T)\mapsto t\iv f(t\iv)\in L^2(\T).
$$
The image of the Riesz projection is equal to the Hardy space
$$
H^2(\T)=\Big\{\;f\in L^2(\T)\;:\;f_k=0\mbox{ for all } k<0\;\Big\}.
$$

For $a\in \Li$ the Toeplitz and Hankel operators are bounded linear operators
defined on $H^2(\T)$ by
\be\label{f.THdef}
T(a)=PM(a)P|_{H^2(\T)},\qquad H(a)=PM(a)JP|_{H^2(\T)}.
\ee
The matrix representation of these operators
with respect to the standard basis $\{t^n\}_{n=0}^\iy$ of $H^2(\T)$
is given by infinite Toeplitz and Hankel matrices,
\be\label{f.THmat}
T(a)\cong(a_{j-k})_{j,k=0}^\iy,\qquad
H(a)\cong(a_{j+k+1})_{j,k=0}^\iy.
\ee
The connection to $n\times n$ Toeplitz and Hankel
matrices is given by
\bqn
P_n T(a)P_n\cong T_n(a),\qquad
P_n H(a) P_n\cong H_n(a),
\eqn
where $P_n$ is the finite rank projection operator
\bqn
P_n: \sum_{k\ge0}f_k t^k \in H^2(\T)\mapsto \sum_{k=0 }^{n-1} f_kt^k\in H^2(\T).
\eqn

Toeplitz and Hankel operators satisfy the following well-known formulas,
\bqn
T(ab)&=& T(a)T(b)+H(a)H(\tb),\label{f.Tab}\\[1ex]
H(ab)&=& T(a)H(b)+H(a)T(\tb),\label{f.Hab}
\eqn
where $\tb(t):=b(t\iv)$.
For $a_-\in \ovl{H^\iy(\T)}$ and $a_+\in H^\iy(\T)$ these formulas
specialize to
\be\label{f.THspez}
T(a_-aa_+)=T(a_-)T(a)T(a_+),\qquad
H(a_-a\ta_+)=T(a_-)H(a)T(a_+),
\ee
where
\bqn
H^\iy(\T) &=& \Big\{\;f\in \Li\;:\; f_k=0\mbox{ for all $k<0$}\;\Big\},\nn\\[1ex]
\overline{H^\iy(\T)} &=& \Big\{\;f\in \Li\;:\; f_k=0\mbox{ for all $k>0$}\;\Big\}.\nn
\eqn
A functions $b$ is called {\em even} if $b=\tb$.

We denote by $\W$ the Wiener algebra, which consists of all
function on $\T$ whose Fourier series is absolutely convergent. Moreover, let
\be
\W_+=\W\cap H^\iy(\T),\qquad
\W_-=\W\cap \ovl{H^\iy(\T)}.
\ee
Functions in $\W_\pm$ can be identified with functions
which are analytic in $\{z\in\C:|z|<1\}$ and $\{z\in\C:|z|>1\}\cup\{\iy\}$,
respectively.
The Riesz projection $P$ is bounded on $\W$ and has the
image $\W_+$.

Given a Banach algebra $B$, we denote by $GB$ the group of all invertible elements in $B$.

A sequence of functions $a_n\in \Li$ is said to converge to $a\in \Li$ in measure
if for each $\eps>0$ the Lebesgue measure of the set
$$
\Big\{\;t\in\T\;:\; |a_n(t)-a(t)|\ge \eps\;\Big\}
$$
converges to zero.

A sequence of bounded linear operators $A_n$
on a Banach space $X$ is said to converge strongly on $X$ to an operator $A$ if
$A_n x\to Ax$ for all $x\in X$.

\begin{lemma}\label{l2.7}
Assume that $a_n\in \Li$ are uniformly bounded and converge to $a\in \Li$
in measure. Then
$$
T(a_n)\to T(a)\quad\mbox{ and } \quad
H(a_n) \to H(a)
$$
strongly on $H^2(\T)$, and the same holds for the adjoints.
\end{lemma}
\begin{proof}
If $a_n$ converges in measure to $a$ and is uniformly bounded, then
$a_n$ also converges to $a$ in the $L^2$-norm. Hence for all
$f\in L^\iy$, we have $a_n f\to af$ in the $L^2$-norm. Using an approximation
argument and the uniform boundedness of $a_n$, it follows
that $M(a_n)\to M(a)$ strongly on $L^2(\T)$.
Hence the corresponding Toeplitz and Hankel operators
converge strongly on $H^2(\T)$, too. Since
$T(a_n)^*=T(a_n^*)$ and $H(a_n)^*=H(\tilde{a}_n^*)$, this
holds also for the adjoints.
\end{proof}

An operator $A$ acting on a Hilbert space $H$ is called a trace class operator if it is compact and if the series constituted by the singular values $s_n(A)$
(i.e., the eigenvalues of $(A^*A)^{1/2}$ taking multiplicities into account)
converges. The norm
\bqn
\|A\|_{1} &=& \sum_{n\ge1} s_n(A)
\eqn
makes the set of all trace class operators into a Banach space, which forms
also a two-sided ideal in the algebra of all bounded linear operators on $H$.
Moreover, the estimates
$\|AB\|_1\le\|A\|_1\|B\|$ and $\|BA\|_1\le\|A\|_1\|B\|$ hold,
where $A$ is a trace class operator and $B$ is a bounded operator with the
operator norm $\|B\|$.

If $A$ is a trace class operator, then the operator trace
$\mathrm{trace}(A)$ and the operator determinant $\det (I+A)$ are
well defined. For more information concerning these concepts we refer to
\cite{GK}.

For $r\in [0,1)$ and $\tau \in \T$ we introduce the following
operators $G_{r,\tau}$ acting on $\Li$:
\be\label{f.G}
G_{r,\tau}:a(t)\mapsto b(t)=a\left(\tau\frac{t+r}{1+rt}\right).
\ee
Figuratively speaking, the function $a$ is first stretched at $\tau$ and squeezed at $-\tau$, and then 
rotated on the unit circle such that $\tau$ moves into $1$.
The inverse operator is given by
\bqn\label{f.Giv}
G_{r,\tau}\iv:a(t)\mapsto b(t)=a\left(\frac{t\tau\iv-r}{1-rt\tau\iv}\right).
\eqn

Given $a\in L^\iy(\R)$ we denote by $M_\R(a)$ the multiplication operator
$$
M_\R(a):f(x)\in L^2(\R)\mapsto a(x)f(x)\in L^2(\R)
$$
and by $W_0(a)$ the convolution operator (or, ``two-sided'' Wiener-Hopf operator)
$$
W_0(a)=\mathcal{F} M_\R(a)\mathcal{F}\iv,
$$
where $\mathcal{F}$ stands for the Fourier transform on $L^2(\R)$.
The usual Wiener-Hopf operator and the ``continuous''
Hankel operator acting on $L^2(\R_+)$ are given by
\bqn\label{f.WH}
W(a)&=&M_\R(\chi_{\R_+})W_0(a)M_\R(\chi_{\R_+})|_{L^2(\R_+)},\\[1ex]
\label{f.HaR}
H_\R(a)&=&M_\R(\chi_{\R_+})W_0(a)\hat{J}M_\R(\chi_{\R_+})|_{L^2(\R_+)},
\eqn
where $(\hat{J}f)(x)=f(-x)$ and $\chi_{\R_+}$ is the characteristic
function of the positive real half axis $\R_+$. If $a\in L^1(\R)$, then
$W(a)$ and $H_\R(a)$ are integal operators on $L^2(\R)$ with the
kernel $\hat{a}(x-y)$ and $\hat{a}(x+y)$, respectively, where
$$
\hat{a}(\xi)=\frac{1}{2\pi}\int_{-\iy}^\iy e^{-ix\xi} a(x)\,dx
$$
stands for the Fourier transform of $a$.

It is important to note that Wiener-Hopf and continous Hankel operators are
related to their discrete analogues by a unitary transform
$S:H^2(\T)\to L^2(\R_+)$,
\bqn\label{f.cd}
W(a)=S T(b) S^*,\qquad
H_{\R}(a)= SH(b)S^*,
\eqn
where
\bqn\label{f.frac}
a(x)=b\left(\frac{1+ix}{1-ix}\right).
\eqn
Let $\Pi_\al$ stand for the projection operator,
\bqn\label{f.Pi}
\Pi_\al&:& f(t)\in L^2(\R_+)\mapsto \chi_{[0,\al]}(x)f(x)\in L^2(\R_+).
\eqn
The image of $\Pi_\al$ can be identified with the space $L^2[0,\al]$.


\section{Outline of the proof}

The main idea of the proof is to establish the identity
\bqn
\det(I-K_{2\al}) &=&
\exp\left(-\frac{\al^2}{2}\right)
\det\Big[\Pi_\al(I+H_\R(\hat{u}_{-1/2}))\iv \Pi_\al\Big]\nn\\
&&\times
\det\Big[\Pi_\al(I-H_\R(\hat{u}_{1/2}))\iv \Pi_\al\Big]
\label{f.id1}
\eqn
and then to apply results of \cite{BE3} in order to determine
the asymptotics of the two operator determinants appearing therein.
Notice that the exponential part already contains the leading
term of the asymptotics of $\det(I-K_{2\al})$.

The appearance of operator determinants
$$
\det\Big[\Pi_\al(I\pm H_\R(\hat{u}_{\beta}))\iv \Pi_\al\Big]
$$
might be quite unmotivated.  Therein, $\hat{u}_\beta\in L^\iy(\R)$ is the function 
defined by
\bqn\label{f.uhat}
\hat{u}_\beta(x) &=& \left(\frac{x-i}{x+i}\right)^\beta,\qquad x\in\R,
\eqn
which is continuous on $\R\setminus \{0\}$, has a jump discontinuity at
$x=0$ and the limits $u(\pm\iy)=1$.
The operators $\Pi_\al(I\pm H_\R(\hat{u}_\beta))\iv \Pi_\al$ are understood as 
operators acting on $L^2[0,\al]$ and then the operator determinants are well-defined.
As is shown in \cite{BE3}, these determinants are related to determinants of Wiener-Hopf-plus-Hankel operators $W(\hat{v}_\beta)\pm H_\R(\hat{v}_\beta)$, where the
symbol $\hat{v}_\beta(x)=(x^2/(1+x^2))^\beta$ has a singularity at $x=0$. The precise relationship is
$$
\det \Pi_\al(W(\hat{v}_\beta)\pm H_\R(\hat{v}_\beta))\Pi_\al = 
e^{-\beta\al}\det\Pi_\al(I\pm H_\R(\hat{u}_{-\beta}))\iv \Pi_\al.
$$
We will not make any further comments, but refer to \cite{BE3} for more information.

The proof of identity (\ref{f.id1}) is accomplished by discretizing the sine kernel operator, which yields a Toeplitz operator, and by  making use of identities between determinants of symmetric Toeplitz matrices, of Hankel matrices and of symmetric Toeplitz-plus-Hankel matrices.
These identities have been established in \cite{BE1}. Unfortunately, these
identities cannot always be applied directly. Thus further comments might
be useful.

The reason why the determinant $\det(I-K_\al)$ is very hard to deal with is that $I-K_\al$ is a finite Wiener-Hopf operator whose generating function vanishes on the whole interval $[-1,1]$. Indeed,  $I-K_\al=\Pi_\al W(1-\chi_{[-1,1]}) \Pi_\al$. Let us remark that the discrete analogue has been 
studied by Widom \cite{W1}. He considered
Toeplitz determinants with a generating function which is even and
vanishes on a single subarc of $\T$, but is elsewhere non-zero and smooth. In the Wiener-Hopf case no comparable result is known so far.

A straightforward discretization of $I-K_\al$ leads to a Toeplitz operator
$T_n(\chi_{\frac{\al}{n}})$ whose (even) generating function $\chi_{\frac{\al}{n}}$ vanishes on a subarc of $\T$, but depends on $n$.
For this reason, the results of \cite{W1} cannot be applied.
This is the place where we use a result of \cite{BE1}. We can identify the
Toeplitz determinant with a determinant of a Hankel operator
$\det H_n[\hat{b}_{\al,n}]$. The crucial point is that although
the function $\hat{b}_{\al,n}$ is not supported on all of $[-1,1]$ it
is supported on a symmetric subinterval. It is thus possible to
pull out a factor of this Hankel determinant (which gives precisly the leading exponential part referred to above after taking the limit $n\to\iy$) to obtain
a Hankel determinant $\det H_n[b_{\al,n}]$ whose generating function
$b_{\al,n}$ is supported on all of $[-1,1]$.

Unfortunately, the function $b_{\al,n}$ is of such a form that one cannot go back to a
Toeplitz determinant by the results of \cite{BE1}. However, another result of
\cite{BE1} establishes an identity between a Hankel determinant
$H_n[b]$ and a determinant of a Toeplitz-plus-Hankel matrix
$\det(T_n(a)+H_n(a))$. Still this result cannot be applied directly.
First of all, the assumptions for this identity are not fulfilled. But worse,
proceeding formally would lead to a function $a$ which does not
belong to $L^1(\T)$.

The way out of this situation is accomplished by establishing the identity
$$
\det(T_n(a)+H_n(a))=\det\Big[P_n (I+H(\psi))\iv P_n\Big]
$$
for nicely behaved functions $a$ and $\psi$.
This allows us to derive the identity
$$
\det H_n[b]=\det\Big[P_n (I+H(\psi))\iv P_n\Big]
$$
for nicely behaved $b$ and $\psi$, and thus we have a possibility of by-passing Toeplitz-plus-Hankel
determinants. We are then able to approximate the generating functions
$b_{\al,n}$ by nicely behaved functions $b$ and obtain
an identity
\be\label{f.id2}
\det H_n[b_{\al,n}]=\det\Big[P_n (I+H(\psi_{\al,n}))\iv P_n\Big].
\ee
Therein the function $\psi_{\al,n}$ is a certain
piecewise continuous function. We remark that the approximation argument is not
quite easy to establish. It requires in particular a stability analysis
for which we resort to results of \cite{ES1}.
The derivation of (\ref{f.id2}) will
be presented in Section \ref{sec4:2}. Before, in Section \ref{sec4:1},
we prove that the operators $I+H(\psi)$ are invertible for just those
functions $\psi$ for which this assertion is needed.

Hence we end up with tackling with the determinant
$$
\det\Big[P_n (I+H(\psi_{\al,n}))\iv P_n\Big].
$$
Analyzing the function $\psi_{\al,n}$ one notices that
$\psi_{\al,n}$ converges uniformly on compact subsets of  $\T\setminus\{-1,1\}$ to a constant
function (for which the Hankel operators would vanish). However, near $t=1$ and $t=-1$, this function shows a considerably more complicated behavior.
Still, one can separate these singularities and prove that the above determinant
behaves asymptotically like
$$
\det\Big[P_n (I+H(\psi_{\al,n}^{(1)})\iv P_n\Big]
\times
\det\Big[P_n (I+H(\psi_{\al,n}^{(-1)})\iv P_n\Big],
$$
where the functions $\psi_{\al,n}^{(1)}$ and $\psi_{\al,n}^{(-1)}$ have ``singular'' behavior only at $t=1$ and $t=-1$, respectively.
Now one can examine these two determinant separately, where in the limit
$n\to\iy$ one arrives at the continuous analogues of these determinant,
$$
\det\Big[\Pi_\al(I+H_\R(\hat{u}_{-1/2}))\iv \Pi_\al\Big]
\quad\mbox{ and }\quad
\det\Big[\Pi_\al(I-H_\R(\hat{u}_{1/2}))\iv \Pi_\al\Big].
$$
The different sign in front of the Hankel operators comes
from the fact that $\psi_{\al,n}^{(1)}$ has its singularity at
$t=1$, while $\psi_{\al,n}^{(-1)}$ has its singularity at
$t=-1$.
The proof of the separation of the singularities as well as
the last step requires a couple of technical result, in which we have to prove
that certain operators converge in trace class norm. These results will be
established in Section \ref{sec4:3}.

The actual proof of the asymptotic formula  as it has been outlined here (i.e., mainly the identity (\ref{f.id1})) will
be given in Section \ref{sec:5}.


\section{Auxiliary results}

\subsection{Invertibility of certain operators $I+H(\psi)$}
\label{sec4:1}

In this section we are going to prove that operators of the
form $I+H(\psi)$ for certain concrete (piecewise continuous)
functions $\psi$ are invertible.

For $\tau\in\T$ and $\beta\in\C$ we introduce the functions
\bqn
\eta_{\beta,\tau}(t)=(1-t/\tau)^{\beta},\qquad
\xi_{\beta,\tau}(t)=(1-\tau/t)^{\beta},
\eqn
where these functions are analytic in an open neighborhood of
$\{\;z\in\C\;:\;|z|\le 1,\;z\neq \tau\}$ and
$\{\;z\in\C\;:\;|z|\ge 1,\;z\neq \tau\}\cup\{\iy\}$, respectively,
and the branch of the power function is chosen in such a way that
$\eta_{\beta,\tau}(0)=1$ and $\xi_{\beta,\tau}(\iy)=1$.
We also need the function
\bqn\label{f.u}
u_{\beta,\tau}(e^{i\theta})=\exp(i\beta(\theta-\theta_0-\pi)),
\quad 0<\theta-\theta_0<2\pi,\quad \tau=e^{i\theta_0},
\eqn
which is continuous on $\T\setminus\{\tau\}$ and has a jump discontinuity
at $t=\tau$.
Notice that
\bqn
u_{\beta,\tau}(t)=\eta_{\beta,\tau}(t)\xi_{-\beta,\tau}(t),\qquad
u_{\beta+n,\tau}(t)=(-t/\tau)^n u_{\beta,\tau}(t).
\eqn

The essential spectrum $\mathrm{sp}_{\mathrm{ess}} A$
of a bounded linear operator $A$ defined on a
Banach space is the set of all $\lambda \in \C$ for which
$A-\lambda I$ is not a Fredholm operator.

We also introduce the Hardy space
\bqn
\ovl{H^2(\T)}=\Big\{\;f\in L^2(\T)\;:\; f_k=0 \mbox{ for all } k>0\;\Big\}.
\eqn
which consists of those functions $f$ for which $\bar{f}\in H^2(\T)$.
Notice that $f\in H^2(\T)$ if and only if $\tilde{f}\in \ovl{H^2(\T)}$.

\begin{proposition}\label{p2.2}
The following operators are invertible on  $H^2(\T)$:
$$
\begin{array}{rclcrcl}
A_1&=&I+H(u_{-1/2,1}),&\qquad&
A_2&=&I-H(u_{1/2,1}),\\[1ex]
A_3&=&I-H(u_{-1/2,-1}),&\qquad&
A_4&=&I+H(u_{1/2,-1}).\\[1ex]
\end{array}
$$
\end{proposition}
\begin{proof}
Let us first consider the operators $A_1$ and $A_2$.
We use a result of Power \cite{Po1} in order to determine the
essential spectrum of a Hankel operator with piecewise continuous symbol.
It says that the essential spectrum is a union of intervals in the complex
plane, namely
\bqn\label{f.Spess}
\mathrm{sp}_{\mathrm{ess}}H(b) &=&[0,ib_{-1}]\cup[0,-ib_1]\cup
\bigcup_{\tau\in\T_+}
\left[-i\sqrt{b_\tau b_{\bar{\tau}}},i \sqrt{b_\tau b_{\bar{\tau}}} \right].
\eqn
Therein we use the notation $b_\tau=(b(\tau+0)-b(\tau-0))/2$ with
$b(\tau\pm0)=\lim_{\eps\to\pm0} b(\tau e^{i\eps})$, and
$\T_+:=\{\tau\in\T\;:\;\Im \tau>0\}$. This result can also be obtained from the
more general results contained in \cite{Po2} and \cite[Secs.~4.95-4.102]{BS}.

Clearly, for our functions $b=u_{-1/2,1}$ and $b=u_{1/2,1}$, respectively,
we have $b_\tau=0$ for $\tau\neq1$. For $\tau=1$, we obtain
$b_1=i$ in case $b=u_{-1/2,1}$ and $b_1=-i$ in case $b=u_{1/2,1}$,
respectively. Hence it follows that
$$
\mathrm{sp}_{\mathrm{ess}}H(u_{-1/2,1})=[0,1]\quad\mbox{ and }\quad
\mathrm{sp}_{\mathrm{ess}}H(u_{1/2,1})=[0,-1].
$$
{}From this we conclude that $I+H(u_{-1/2,1})$ and $I-H(u_{1/2,1})$ are Fredholm
operators with Fredholm index zero.

In order to prove invertibility it thus suffices to show that the kernels of these
operators are trivial.

Let $f_+\in H^2(\T)$ belong to the kernel of $I+H(u_{-1/2,1})$. Then
$$
f_+(t)+u_{-1/2,1}(t)t\iv\tilde{f}_+(t) =:f_-(t)\in t\iv \overline{H^2(\T)}.
$$
Using the identity
$u_{-1/2,1}(t)=-t\iv u_{1/2,1}(t)=-t\iv\eta_{1/2,1}(t)\xi_{-1/2,1}(t)$, we obtain
$$
f_0(t):= t\xi_{1/2,1}(t)f_+(t)-t\iv \eta_{1/2,1}(t)\tilde{f}_+(t)
=t\xi_{1/2,1}(t) f_-(t).
$$
{}From the definition of $f_0$ it follows that $f_0=-\tilde{f}_0$, while
$t\xi_{1/2,1}(t) f_-(t)\in \overline{H^2(\T)}$.
Hence $f_0=0$ and we have shown that
$$
t\xi_{1/2,1}(t)f_+(t)-t\iv \eta_{1/2,1}(t)\tilde{f}_+(t)=0.
$$
This implies
$$
f_+(t)=u_{1/2,1}(t) t^{-2} \tilde{f}_+(t)=-u_{-1/2,1}(t) t\iv \tilde{f}_+(t).
$$
Using $u_{-1/2,1}(t)=\eta_{-1/2,1}(t)\xi_{1/2,1}(t)$ we obtain
$$
\eta_{1/2,1}(t)f_+(t)=-\xi_{1/2,1}(t) t\iv\tilde{f}_+(t).
$$
Therein the left hand side belongs to $H^2(\T)$ whereas the right hand side belongs
to $t\iv \overline{H^2(\T)}$. Hence they must be zero. This implies $f_+=0$
as desired.

Now let $f_+\in H^2(\T)$ belong to the kernel of $I-H(u_{1/2,1})$. Then
$$
f_+(t)-u_{1/2,1}(t)t\iv \tilde{f}_+(t)=:f_-(t)\in t\iv \overline{H^2(\T)}.
$$
Using $u_{1/2,1}(t)=\eta_{1/2,1}(t)\xi_{-1/2,1}(t)$ we obtain
$$
f_0(t):=\xi_{1/2,1}(t)f_+(t)-\eta_{1/2,1}(t) t\iv \tilde{f}_+(t)=\xi_{1/2,1}(t)f_-(t).
$$
Since $f_0(t)=-t\iv \tilde{f}_0(t)$ and
$\xi_{1/2,1}(t)f_-(t)\in t\iv \overline{H^2(\T)}$, we immediately obtain $f_0=0$.
Hence
$$
f_+(t)=u_{1/2,1}(t) t \iv \tilde{f}_+(t).
$$
Using the formula
$u_{1/2,1}(t) t\iv =-u_{-1/2,1}(t)=-\eta_{-1/2,1}(t)\xi_{1/2,1}(t)$,
we conclude that
$$
\eta_{1/2,1}(t) f_+(t)=\xi_{1/2,1}(t)\tilde{f}_+(t).
$$
The left hand side belongs to $H^2(\T)$, whereas the right hand side belongs to
$\ovl{H^2(\T)}$. Hence it is a constant which must be zero because
$\xi_{-1/2,1}\notin L^2(\T)$. Thus we obtain $f_+=0$, which proves that the kernel
is trivial.

Finally, we can say that the operators $A_3$ and $A_4$ can be treated analogously.
However, we can also conclude their invertibility directly by remarking that
$A_3=WA_1 W$ and $A_4=WA_2W$, where $(Wf)(t)=f(-t)$, $t\in\T$.
\end{proof}

Next we introduce the function
\bqn\label{f.chi.def}
\chi(e^{i\theta})= \left\{\ba{rl} i &\mbox{ if } 0<\theta<\pi\\[1ex]
-i&\mbox{ if } -\pi<\theta<0.\ea\right.
\eqn
It is easy to see that the identity
\be\label{f.chi}
\chi(t)=u_{-1/2,1}(t) u_{1/2,-1}(t)=-u_{1/2,1}(t)u_{-1/2,-1}(t)
\ee
holds.

\begin{proposition}\label{p2.4x}
Let $c_+\in G\W_+$ and $\psi(t)=\tc_+(t)c_+\iv(t)\chi(t)$. Then
the operator $I+H(\psi)$ is invertible on $H^2(\T)$.
\end{proposition}
\begin{proof}
The proof goes like the proof of Proposition \ref{p2.2}. First of all we determine
the essential spectrum of $H(\psi)$. Using the notation
$\psi_\tau=(\psi(\tau+0)-\psi(\tau-0))/2$, it easily follows that
$\psi_1=i$, $\psi_{-1}=-i$ and $\psi_\tau=0$ for $\tau\in \T\setminus\{1,-1\}$. Hence by (\ref{f.Spess})
$$
\mathrm{sp}_{\mathrm{ess}} H(\psi)=[0,1],
$$
which implies that $I+H(\psi)$ is a Fredholm operator with index zero. It remains to prove
that the kernel of $I+H(\psi)$ is trivial.

Suppose that $f_+\in H^2(\T)$ belongs to this kernel. Then
$$
f_+(t)+\psi(t)t\iv \tilde{f}_+(t)=:f_-(t)\in t\iv \ovl{H^2(\T)}.
$$
By (\ref{f.chi}) we can write
$$
\chi(t)=-t\iv u_{1/2,1}(t) u_{1/2,-1}(t)=-t\iv \xi_{-1/2,1}(t)\xi_{-1/2,-1}(t)
\eta_{1/2,1}(t)\eta_{1/2,-1}(t),
$$
and hence we obtain
\bqn
f_0(t) &:=&t\tc_+\iv(t) \xi_{1/2,1}(t)\xi_{1/2,-1}(t)f_+(t)-
t\iv c_+\iv(t) \eta_{1/2,1}(t)\eta_{1/2,-1}(t)\tilde{f}_+(t)\nn\\[.5ex]
&=&t\tc_+\iv(t)\xi_{1/2,1}(t)\xi_{1/2,-1}(t)f_-(t).\nn
\eqn
Since $f_0=-\tilde{f}_0$ and since the right hand side belongs to
$\ovl{H^2(\T)}$, it follows that $f_0=0$. Thus
$$
f_+(t)+\psi(t)t\iv \tilde{f}_+(t)=0.
$$
Now we write
$$
\chi(t)=t u_{-1/2,1}(t)u_{-1/2,-1}(t)= t \xi_{1/2,1}(t)\xi_{1/2,-1}(t)
\eta_{-1/2,1}(t)\eta_{-1/2,-1}(t),
$$
and it follows
$$
c_+(t)\eta_{1/2,1}(t)\eta_{1/2,-1}(t)f_+(t)=-\tc_{+}(t)
\xi_{1/2,1}(t)\xi_{1/2,-1}(t)\tilde{f}_+(t).
$$
Therein the left hand side belongs to $\ovl{H^2(\T)}$ whereas the right hand
side belongs to $H^2(\T)$. It follows that this expression is zero.
Hence $f_+=0$.
This proves that the kernel is trivial.
\end{proof}


\subsection{A formula for Hankel determinants}
\label{sec4:2}

The goal of this section is to prove the formula
$$
\det H_n[b]= G^n\det\left[P_n\Big(I+H(\psi)\Big)\iv P_n\right]
$$
where $b\in L^1[-1,1]$ is a (sufficiently smooth) continuous nonvanishing function and $\psi$ is a function defined in terms of $b$
(see Theorem \ref{t2.4} below).

In the following proposition, we denote by $\log a$ any continuous logarithm
of the function $a$ and by $[\log a]_n$ the $n$-th Fourier coefficient of
$\log a$.

\begin{proposition}\label{p2.1}
Let $a\in G\W$ be an even function.
Then there exists a function $a_+\in G\W_+$ with $a_+(0)=1$ such that
\bqn\label{f.fact}
a(t)&=& \ta_+(t)G a_+(t), \qquad t\in \T,
\eqn
where $G=\exp([\log a]_0)$. Moreover, the operator $I+H(\psi)$
is invertible on $H^2(\T)$, where
\bqn
\psi(t) &=& \ta_+(t)a_+\iv(t),
\eqn
and for all $n\ge1$ the following identity holds:
\bqn\label{f.BO}
\det\Big(T_n(a)+H_n(a)\Big) &=& G^n \det\left[ P_n\Big(I+H(\psi)\Big)\iv P_n\right]
\eqn
\end{proposition}
\begin{proof}
An even continuous nonvanishing function has winding number zero and thus
possesses a continuous logarithm. Since $a\in \W$ it is easy to see
(e.g., by an approximation argument and by using the fact that
$\W$ is a Banach algebra in which the trigonometric polynomials are dense)
that $\log a \in \W$.
We define
\bqn
(\log a)_+(t) &:=& \sum_{n=1}^\iy [\log a]_n t^{n},\nn
\eqn
which belongs to $\W_+$ since $\W_+$ is the image
of the Riesz projection $P$ on the space $\W$.
Note that $(\log a)(t)=(\log a)_+(t\iv)+[\log a]_0+(\log a)_+(t)$
since $\log a$ is also an even function. Upon defining
$a_+=\exp(\log a)_+$, which belongs to the Banach algebra $\W_+$,
the factorization (\ref{f.fact}) follows immediately.

Now we employ formulas (\ref{f.Tab}) and (\ref{f.Hab}) in connection with
$a=\ta$  to conclude that
$$
\Big(T(a)+H(a)\Big)\Big(T(a\iv)+H(a\iv)\Big)=
\Big(T(a\iv)+H(a\iv)\Big)\Big(T(a)+H(a)\Big)=I.
$$
Moreover, using formulas  (\ref{f.THspez}) we deduce that
$$
T(a\iv)+H(a\iv)=T(\ta_+\iv G\iv)\Big(I+H(\psi)\Big)T(a_+\iv).
$$
The just proved invertibility of $T(a\iv)+H(a\iv)$ implies that
$I+H(\psi)$ is invertible since $T(\ta_+\iv)$ and $T(a_+\iv)$ are
invertible by (\ref{f.THspez}). It follows that
$$
T(a)+H(a)=\Big(T(a\iv)+H(a\iv)\Big)\iv=
T(a_+)\Big(I+H(\psi)\Big)\iv T(\ta_+G),
$$
whence we obtain
$$
T_n(a)+H_n(a)=T_n(a_+)\left[P_n\Big(I+H(\psi)\Big)\iv P_n\right] T_n(\ta_+G)
$$
since $T(a_+)$ and $T(\ta_+)$ are lower and upper, respectively, triangular matrices
in the standard matrix representation. Noting that
the diagonal entries of $T_n(a_+)$ and $T_n(\ta_+)$ are equal
to $[a_+]_0=a_+(0)=1$ implies assertion (\ref{f.BO}) by taking the determinant.
\end{proof}

Next we cite the following result from \cite[Thm.~2.3]{BE1}.
Recall the definition of the Hankel operator $H_n[b]$ given in (\ref{f.Hn2}).

\begin{proposition} \label{p2.3}
Let $a\in L^1(\T)$ be an even function and let $b\in L^1[-1,1]$ be given by
\bqn\label{f.TH-H}
b(\cos\theta)=a(e^{i\theta})\sqrt{\frac{1+\cos\theta}{1-\cos\theta}}.
\eqn
Then $\det\Big(T_n(a)+H_n(a)\Big)=\det H_n[b]$.
\end{proposition}

We remark in this connection that under the assumption (\ref{f.TH-H})
we have $b\in L^1[-1,1]$ if and only if
$a(e^{i\theta}) (1+\cos\theta)\in L^1(\T)$.

In regard to the following theorem recall the definition (\ref{f.chi.def}) of the function $\chi$. Moreover, notice that
$I+H(\psi)$ is invertible by Proposition \ref{p2.4x}.

\begin{theorem}\label{t2.4}
Let $c\in G\W$ be an even function such that
\bqn
c(t)&=& \tc_+(t)G c_+(t)
\eqn
where $c_+\in G\W_+$, $c_+(0)=1$, and $G=\exp([\log c]_0)$.
Moreover, let
$$
b(\cos\theta)=c(e^{i\theta}),\qquad
\psi(e^{i\theta})=\tc_+(e^{i\theta}) c_+\iv(e^{i\theta}) \chi(e^{i\theta}).
$$
Then
$$
\det H_n[b]= G^n\det\left[P_n\Big(I+H(\psi)\Big)\iv P_n\right].
$$
\end{theorem}
\begin{proof}
The proof will be carried out by an approximation argument. For $r\in[0,1)$
introduce the even function
$$
a_r(t)= c(t)\sqrt{\frac{(1-rt)(1-rt\iv)}{(1+rt)(1+rt\iv)}},\qquad t\in\T.
$$
The function corresponding to $a_r$ by means of (\ref{f.TH-H}) is then
$$
b_r(x)=b(x)\sqrt{\frac{2+2x}{1+r^2+2rx}}\sqrt{\frac{1+r^2-2rx}{2-2x}},\qquad
x\in(-1,1).
$$
It is easy to verify that $b_r\to b$ in the norm of $L^1[-1,1]$.
Hence (for each fixed $n$)
$$
H_n[b]=\lim_{r\to1} H_n[b_r]=\lim_{r\to1}\det\Big(T_n(a_r)+H_n(a_r)\Big)
$$
by Proposition \ref{p2.3}.

The canonical Wiener-Hopf factorization of $a_r$ is given
by $a_r(t)=\ta_{r,+}(t) G a_{r,+}(t)$
with
$$
a_{r,+}(t)=c_+(t)\frac{(1-rt)^{1/2}}{(1+rt)^{1/2}}.
$$
Upon putting
$$
\psi_r(t)=\ta_{r,+}(t)a_{r,+}\iv(t)=
\tc_{+}(t)c_{+}\iv(t)\left(\frac{1-rt}{1-rt\iv}\right)^{-1/2}
\left(\frac{1+rt}{1+rt\iv}\right)^{1/2},
$$
we conclude from Proposition \ref{p2.1} that
$$
\det\Big(T_n(a_r)+H_n(a_r)\Big)=G^n\det\left[P_n \Big( I+H(\psi_r)\Big)\iv P_n\right].
$$
Hence
$$
\det H_n[b]= G^n\lim_{r\to 1} \det\left[P_n \Big( I+H(\psi_r)\Big)\iv P_n\right].
$$
Since
\be\label{f.fpm}
f_r^\pm := \left(\frac{1\mp rt}{1\mp rt\iv}\right)^{\mp 1/2}\to u_{\mp 1/2,\pm 1}(t)
\ee
in measure as $r\to 1$ and because of (\ref{f.chi}),
it follows that $\psi_r\to \psi$ in measure. Since the sequence is bounded in
the $L^\iy$-norm, it follows that $H(\psi_r)$ converges strongly to $H(\psi)$
on $H^2(\T)$ by Lemma \ref{l2.7}.

In order to conclude that
\bqn\label{f.sc}
(I+H(\psi_r))\iv & \to &(I+H(\psi))\iv
\eqn
strongly on $H^2(\T)$,
it is necessary and sufficient that (for some $r_0\in[0,1)$)
$$
\sup_{r\in[r_0,1)}\left\|(I+H(\psi_r))\iv\right\| < \iy.
$$

In order to analyse this stability condition we apply the results of
\cite[Secs.~4.1--4.2]{ES1}.  These results establish the existence of certain
mappings $\Phi_0$ and $\Phi_\tau$, $\tau\in\T$, which in our case evaluate
as follows. Recall the definition (\ref{f.G}) of the mapping $G_{r,\tau}$.
Because of (\ref{f.fpm}) we have
$$
\Phi_0[f_r^\pm]:=\mu\mbox{-}\lim_{r\to1} f_r^\pm =u_{\mp1/2,\pm1}
$$
where $\mu\mbox{-}\lim$ stands for the limit in measure. Furthermore since
$f_r^\pm \to u_{\mp1/2,\pm1}$ locally uniformly on $\T\setminus\{\pm1\}$, we have
$$
\Phi_\tau[f_r^\pm]= \mu\mbox{-}\lim_{r\to1} G_{r,\tau} f_r^\pm = u_{\mp 1/2,\pm 1}(\tau)
$$
for $\tau\neq\pm 1$. Finally,
$$
\Phi_{\pm 1}[f_r^\pm]= \mu\mbox{-}\lim_{r\to 1} G_{r,\pm 1} f_r^\pm=
\mu\mbox{-}\lim_{r\to 1}
\left(\frac{1+rt}{1+rt\iv}\right)^{\pm1/2} = u_{\pm 1/2,-1}.
$$
Since $\psi_r =\tc_+ c_+\iv f_r^+f_r^-$ we conclude
\bqn
\Phi_0[\psi_r] &=& \tc_+c_+\iv u_{-1/2,1}u_{1/2,-1}=\psi,\nn\\[1ex]
\Phi_1[\psi_r] &=& u_{1/2,-1},\nn\\[1ex]
\Phi_{-1}[\psi_r] &=& u_{-1/2,-1},\nn\\[1ex]
\Phi_\tau[\psi_r] &=& \mbox{constant function, }\quad\tau\in \T\setminus\{-1,1\}.\nn
\eqn
The stability criterion in \cite{ES1} (Thm.~4.2 and Thm.~4.3) says that $I+H(\psi_r)$ 
is stable if and only if the operators
\begin{itemize}
\item[(i)] $\Psi_0[I+H(\psi_r)]=I+H(\Phi_0[\psi_r])=I+H(\psi)$,
\item[(ii)] $\Psi_1[I+H(\psi_r)]=I+H(\Phi_1[\psi_r])=I+H(u_{1/2,-1})$,
\item[(iii)] $\Psi_{-1}[I+H(\psi_r)]=I-H(\Phi_{-1}[\psi_r])=I-H(u_{-1/2,-1})$,
\item[(iv)] $\Psi_\tau[I+H(\psi_r)]=$
$$
\twomat{I&0\\0&I}+\twomat{P&0\\0&Q}\twomat{M(\Phi_\tau[\psi_r])&0\\0&
M(\widetilde{\Phi_{\bar{\tau}}[\psi_r]})}\twomat{0&I\\ I&0}
\twomat{P&0\\0&Q}=\twomat{I&0\\0&I}
$$
($\tau\in\T$, $\mathrm{Im}(\tau)>0$)
\end{itemize}
are invertible. Clearly, by Proposition \ref{p2.2} and Proposition \ref{p2.4x} this is the case.
Hence the sequence $I+H(\psi_r)$ is stable and (\ref{f.sc}) follows.
We conclude (for fixed $n$) that the $n\times n$ matrices
$P_n(I+H(\psi_r))\iv P_n$ converge to $P_n (I+H(\psi))\iv P_n$
as $r\to 1$,
whence it follows that the corresponding  determinants converge, too.
This completes the proof.
\end{proof}


\subsection{Convergence in trace class norm}
\label{sec4:3}

In this section we are going to prove a couple of technical results.
We are mainly concerned with proving that certain sequences converge in the trace norm.

Let $\pcea$ stand for the set of all functions on $\T$ which are
absolutely continuous on $\T\setminus\{-1,1\}$ and which possess
one-sided limits at $t=1$ and $t=-1$.

\begin{lemma}\label{l2.6}
Let $a\in C(\T)$ be a function such that $a'\in\pcea$.
Then $H(a)$ is a trace class operator on $H^2(\T)$ and
\bqn\label{f.31y}
\|H(a)\|_1 &\le & C\Big(\|a\|_{L^\iy(\T)}+\|a'\|_{L^\iy(\T)}+\|a''\|_{L^1(\T)}\Big).
\eqn
\end{lemma}
\begin{proof}
{}From partial integration it follows that the Fourier coefficients
$a_k$ are $O(k^{-2})$ as $k\to\iy$, where the constant involved in this estimate
is given in terms of the norms of $a$, $a'$ and $a''$. We write the operator $H(a)$
as a product $AB$ with operators $A$ and $B$ given by its matrix representation with respect
to the standard basis by
$$
A=\Big(a_{j+k+1}(1+k)^{1/2+\eps}\Big)_{j,k=0}^\iy,\quad
B=\diag\Big((1+k)^{-1/2-\eps}\Big)_{k=0}^\iy.
$$
Both $A$ and $B$ are Hilbert-Schmidt operators if $0<\eps<1/2$ with straightforward
estimates for their norms. Hence $H(a)$ is a trace class operator, whose
norm can be estimated by (\ref{f.31y}).
\end{proof}

In the following proposition we prove that certain operators converge to zero
in the trace norm. The proof is very technical. It might be illustrative to remark that the convergence of these operators in the operator norm is
almost obvious.

Recall the definition of the operator $G_{\mu,\tau}\iv$ given in (\ref{f.Giv}).

\begin{proposition}\label{p2.8}
Let
\be\label{phi.mu}
\psi_\mu^{(1)} = G_{\mu,1}^{-1}(u_{-1/2,1}-1),\qquad
\psi_\mu^{(-1)} = G_{\mu,-1}^{-1}(u_{1/2,1}-1)
\ee
with $\mu\in[0,1)$. Then
the operators
$$
H(\psi_{\mu}^{(1)})H(\psi_{\mu}^{(-1)}),
\qquad
H(\psi_{\mu}^{(-1)})H(\psi_{\mu}^{(1)}),
\quad\mbox{ and }\quad
H(\psi_{\mu}^{(1)}\psi_{\mu}^{(-1)})
$$
are trace class operators and converge to zero in the trace norm
as $\mu\to1$.
\end{proposition}
\begin{proof}
Let us first notice that (with the proper choice of the square-root),
\be\label{f.26}
\psi_{\mu}^{(1)}(t)=\left(-\frac{t-\mu}{1-\mu t}\right)^{-1/2}-1,\qquad
\psi_{\mu}^{(-1)}(t)=\left(\frac{t+\mu}{1+\mu t}\right)^{1/2}-1.
\ee
In particular, $\psi_{\mu}^{(1)}$ has a jump discontinuity at $t=1$ and
vanishes at $t=-1$ while $\psi_{\mu}^{(-1)}$ has a jump discontinuity at $t=-1$
and vanishes at $t=1$. Moreover, both functions are uniformly bounded and
\be\label{f.27}
\psi_{\mu}^{(1)}\to 0,\qquad
\psi_{\mu}^{(-1)}\to 0,
\ee
uniformly on each compact subset of $\T\setminus\{1\}$ and $\T\setminus\{-1\}$, respectively.

In order to prove the assertion for the operator
$H(\psi_{\mu}^{(1)})H(\psi_{\mu}^{(-1)})$, let $f$ and $g$
be smooth functions on $\T$ with $f+g=1$ such that
$f(t)$ vanishes identically in a neighborhood of $1$
(say for $|\arg t|\le \pi/3$)
and $g(t)$ vanishes identically in a neighborhood of $-1$
(say for $|\arg t|\ge  2\pi/3$). Then (see (\ref{f.Hab}))
\bqn
H(\psi_{\mu}^{(1)})H(\psi_{\mu}^{(-1)}) &=&
H(\psi_{\mu}^{(1)})T(f)H(\psi_{\mu}^{(-1)}) +
H(\psi_{\mu}^{(1)})T(g)H(\psi_{\mu}^{(-1)})\nn\\[.5ex]
&=&
H(\psi_{\mu}^{(1)}\tilde{f})H(\psi_{\mu}^{(-1)}) -
T(\psi_{\mu}^{(1)})H(\tilde{f})H(\psi_{\mu}^{(-1)}) \nn\\
&&+
H(\psi_{\mu}^{(1)})H(g\psi_{\mu}^{(-1)})-
H(\psi_{\mu}^{(1)})H(g)T(\widetilde{\psi_{\mu}^{(-1)}}).\nn
\eqn
Clearly, $H(\tilde{f})$ and $H(g)$ are trace class operators. Due to the afore-mentioned
fact that $\psi_\mu^{(1)}$ and $\psi_\mu^{(-1)}$ are uniformly bounded and
because of the convergence (\ref{f.27}), Lemma \ref{l2.7}
implies that the operators
$$
H(\psi_\mu^{(1)}),\quad T(\psi_\mu^{(1)}),\quad
H(\psi_\mu^{(-1)}),\quad T(\widetilde{\psi_\mu^{(1)}})
$$
and their adjoints converge strongly to zero as $\mu\to1$.
We can conclude that $H(\psi_\mu^{1})H(\psi_\mu^{(-1)})$
is a trace class operator and converges in the trace norm to zero as
soon as we have shown that
$$
H(\psi_\mu^{(1)}\tilde{f})\quad\mbox{ and }\quad H(g\psi_\mu^{(-1)})
$$
are trace class operators which converge to zero in the trace norm.
On account of Lemma \ref{l2.6} this is true if
$$
\psi_\mu^{(1)}\tilde{f}\in C(\T),\quad
(\psi_\mu^{(1)}\tilde{f})'\in \pcea,
$$
if
$$
\|\psi_\mu^{(1)}\tilde{f}\|_{L^\iy}\to0,\quad\|(\psi_\mu^{(1)}\tilde{f})'\|_{L^\iy}\to0,\quad
\|(\psi_\mu^{(1)}\tilde{f})''\|_{L^1}\to0
$$
and if similar statements hold for $g\psi_\mu^{(-1)}$.
Due to the fact that $f$ vanishes on a neighborhood of $1$, these conditions
are fulfilled if
\be\label{f.28}
\psi_\mu^{(1)}|_{\Ta}\in C(\Ta),\quad
(\psi_\mu^{(1)})'|_{\Ta}\in C(\Ta),\quad
(\psi_\mu^{(1)})''|_{\Ta}\in C(\Ta),
\ee
and if
\be\label{f.29}
\|\psi_\mu^{(1)}|_{\Ta}\|_{L^\iy(\Ta)}\to 0,\quad
\|(\psi_\mu^{(1)})'|_{\Ta}\|_{L^\iy(\Ta)}\to 0,\quad
\|(\psi_\mu^{(1)})''|_{\Ta}\|_{L^1(\Ta)}\to 0.
\ee
Therein we have restricted the function $\psi_\mu^{(1)}$ to the interval
$\Ta:=\{t\in\T\;:\;|\arg t|\ge \pi/4\}$.
The corresponding (sufficient) conditions for the function $\psi_\mu^{(-1)}$
are
\be\label{f.30}
\psi_\mu^{(-1)}|_{\Tb}\in C(\Tb),\quad
(\psi_\mu^{(-1)})'|_{\Tb}\in C(\Tb),\quad
(\psi_\mu^{(-1)})''|_{\Tb}\in C(\Tb),
\ee
and
\be\label{f.31}
\|\psi_\mu^{(-1)}|_{\Tb}\|_{L^\iy(\Tb)}\to 0,\quad
\|(\psi_\mu^{(-1)})'|_{\Tb}\|_{L^\iy(\Tb)}\to 0,\quad
\|(\psi_\mu^{(-1)})''|_{\Tb}\|_{L^1(\Tb)}\to 0,
\ee
where $\Tb:=\{t\in\T\;:\;|\arg t|\le 3\pi/4\}$.
It is easy to see that conditions (\ref{f.28}) and (\ref{f.30}) and also the
first condition in (\ref{f.29}) and in (\ref{f.31}) are fulfilled.

We will prove the remaining conditions in a few moments, but first we will turn to the convergence of the operators
$H(\psi_\mu^{(-1)})H(\psi_\mu^{(1)})$ and
$H(\psi_\mu^{(1)}\psi_\mu^{(-1)})$. In regard to the operator
$H(\psi_\mu^{(-1)})H(\psi_\mu^{(1)})$ we can proceed analogously and it
turns out that we arrive at the same sufficient conditions
(\ref{f.28})--(\ref{f.31}).

As to the operator $H(\psi_\mu^{(1)}\psi_\mu^{(-1)})$ we have to show
(on account of Lemma \ref{l2.6}) that
\be\label{f.32}
\psi_{\mu}^{(1)}\psi_{\mu}^{(-1)}\in C(\T),\quad
(\psi_{\mu}^{(1)}\psi_{\mu}^{(-1)})'\in \pcea
\ee
and that
\be\label{f.33}
\|\psi_{\mu}^{(1)}\psi_{\mu}^{(-1)}\|_{L^\iy}\to 0,\quad
\|(\psi_{\mu}^{(1)}\psi_{\mu}^{(-1)})'\|_{L^\iy}\to 0,\quad
\|(\psi_{\mu}^{(1)}\psi_{\mu}^{(-1)})''\|_{L^1}\to 0,
\ee
{}From the facts stated at the beginning of the proof it follows that
$\psi_{\mu}^{(1)}\psi_{\mu}^{(-1)}$ is continuous on $\T$ and that
$\psi_{\mu}^{(1)}\psi_{\mu}^{(-1)}$ converges uniformly to zero on $\T$.
Moreover, since the functions $\psi_\mu^{(\pm1)}$ and their derivatives
belong to $\pcea$,
it follows that the derivative of
$\psi_{\mu}^{(1)}\psi_{\mu}^{(-1)}$ of  belongs to $\pcea$, too.
Thus we are left with the proof of the second and third condition in (\ref{f.33}).
We will prove these assertions by separating the singularities at $t=1$ and $t=-1$:
\be
\|(\psi_{\mu}^{(1)}\psi_{\mu}^{(-1)})'|_{\Ta}\|_{L^\iy(\Ta)}\to 0,\quad
\|(\psi_{\mu}^{(1)}\psi_{\mu}^{(-1)})''|_{\Ta}\|_{L^1(\Ta)}\to 0,
\ee
and
\be\label{f.35}
\|(\psi_{\mu}^{(1)}\psi_{\mu}^{(-1)})'|_{\Tb}\|_{L^\iy(\Tb)}\to 0,\quad
\|(\psi_{\mu}^{(1)}\psi_{\mu}^{(-1)})''|_{\Tb}\|_{L^1(\Tb)}\to 0.
\ee

Now turning back to the proof of the yet outstanding conditions in (\ref{f.29})
and (\ref{f.31}), we remark that the interval $\Ta$ can be transformed into
the interval $\Tb$ by a rotation $t\mapsto -t$. This will not precisely transform
the function $\psi_\mu^{(1)}$ into the function $\psi_\mu^{(-1)}$, but into a
similar function of the form (\ref{f.26}), where only the power $1/2$
is replaced by $-1/2$. Without loss of generality we can thus confine ourselves
to the proof of the conditions involving the interval $\Tb$, since the conditions
involving the interval $\Ta$ can be reduced to an analogous situation and can be proved in the same way.

In order to prove (\ref{f.35}) and the last two conditions in (\ref{f.31}) we
use the linear fractional transformation
$$
\sigma(x)=\frac{1+ix}{1-ix},
$$
which maps the extented real line onto the unit circle. Clearly, $\Tb$ corresponds
to
$\sigma\iv(\Tb)=[-1-\sqrt{2},1+\sqrt{2}]=:I_0$. We transform the functions
into
$$
v_\eps(x)=\psi_\mu^{(1)}(\sigma(x)),\qquad
w_\eps(x)=\psi_\mu^{(-1)}(\sigma(x)),
$$
and we also change the parameter $\mu\in[0,1)$ into
$\eps=\frac{1-\mu}{1+\mu}\in(0,1]$. The conditions which we  have to prove
are then equivalent to
\be\label{f.36x}
\|(w_\eps )'|_{I_0}\|_{L^\iy(I_0)}\to0,\qquad
\|(w_\eps )''|_{I_0}\|_{L^1(I_0)}\to0
\ee
and
\be\label{f.37}
\|(v_\eps w_\eps )'|_{I_0}\|_{L^\iy(I_0)}\to0,\qquad
\|(v_\eps w_\eps )''|_{I_0}\|_{L^1(I_0)}\to0
\ee
as $\eps\to0$. Introduce the functions
$$
v(x)=\left(\frac{x-i}{x+i}\right)^{-1/2}-1,\qquad
w(x)=\left(\frac{1+ix}{1-ix}\right)^{1/2}-1,
$$
where $v$ has a jump at $x=0$ and the square-root is chosen such that
$v(\pm\iy)=0$. The function $w$ is continuous on $\R$ with $w(0)=0$ and limits
at $x\to\pm\iy$.
A straightforward computation implies that $v_\eps(x)=v(x/\eps)$
and $w_\eps(x)=w(x\eps)$.

The functions $v$ and $w$ and all of their derivatives are bounded on $\R$.
Thus the conditions in (\ref{f.36x}) follow easily. The function $w$ can be written as
$w(x)=x\tilde{w}(x)$, where $\tilde{w}$ is a function which is locally bounded.
We write
$$
(v_\eps w_\eps)'=v'(x/\eps)x\tilde{w}(\eps x)+v(x/\eps)\eps w'(\eps x)
$$
and see immediately that the second term goes uniformly to zero.
Moreover, $v'(x)\to 0$ as $|x|\to\iy$. Hence $xv'(x/\eps)$
converges uniformly on $I_0$ to zero, which implies that
the first term converges uniformly on $I_0$ to zero.
Thus we have proved that $(v_\eps w_\eps)'$ converges uniformly on $I_0$
to zero as $\eps\to 0$.

Finally, we write the second derivative as
\bqn
(v_\eps w_\eps)'' &=&
\eps^{-1}v''(x/\eps)x\tilde{w}(\eps x)+2v'(x/\eps) w'(\eps x)+
\eps^2 v(x/\eps) w''(\eps x).
\eqn
The $L^1(I_0)$-norm of the first term can be estimated by a constant times
$$
\int_{I_0}|v''(x/\eps) x/\eps|dx \le \eps\int_{\R} |x v''(x)|\, dx,
$$
which converges to zero. The $L^1(I_0)$-norm of the second term can be estimated
by a constant times
$$
\int_{I_0}|v'(x/\eps) |dx \le \eps\int_{\R} | v'(x)|\, dx
$$
and also converges to zero. The last term converges to zero even uniformly.
Hence we have proved the conditions (\ref{f.37}) and the proof is complete.
\end{proof}

In addition to the operators $G_{\mu,\tau}$ we introduce operators
\bqn
R_{\mu,\tau} &:&   f\in H^2(\T)\mapsto g(t)=\frac{\sqrt{1-\mu^2}}{1+\mu t}
G_{\mu,\tau}(f)\in H^2(\T),
\eqn
where $\mu\in[0,1)$ and $\tau\in \{-1,1\}$.

\begin{lemma}\label{l2.9}
For each $\tau\in\{-1,1\}$, the operator $R_{\mu,\tau}$ is unitary on $H^2(\T)$.
Moreover, $R_{\mu,\tau}H(a)R_{\mu,\tau}^*=\tau H(G_{\mu,\tau} a)$ for all
$a\in  \Li$.
\end{lemma}
\begin{proof}
We can define the operators $R_{\mu,\tau}$ also on $L^2(\T)$.
In \cite[Sect.~5.1]{ES1} it is proved that $R_{\mu,\tau}$ are unitary
on $L^2(\T)$ and that
$$
R_{\mu,\tau}PR_{\mu,\tau}^*=P,\quad
R_{\mu,\tau}M(a)R_{\mu,\tau}^*=M(G_{\mu,\tau}\iv a),\quad
R_{\mu,\tau}JR_{\mu,\tau}^*=\tau J.
$$
These statements imply the desired assertions.
\end{proof}

In connection with the following proposition recall that
the operators
$I+H(u_{-1/2,1})$ and $I-H(u_{1/2,1})$ are invertible on $H^2(\T)$
(see Proposition \ref{p2.2}).

Moreover, define the functions
\bqn
\label{f.hal}
h_{\alpha}(t) &=& \exp\left(-\frac{\al(1-t)}{2(1+t)}\right),
\\[1ex]
\label{f.haln}
h_{\al,n}(t)  &=& \left(\frac{t+\mu_{\al,n}}{1+\mu_{\al,n}t}\right)^n,
\eqn
where $\mu_{\al,n}\in [0,1)$ is a sequence for which
\bqn\label{f.mu.an1}
\mu_{\al,n}&=&1-\frac{\al}{2n}+O(n^{-2}),\quad
\mbox{ as } n\to\iy,
\eqn
for each $\al>0$.
Finally, introduce the functions
\be\label{f.psi.pm1}
\psi_{\al,n}^{(1)}= G_{\mu_{\al,n},1}\iv (u_{-1/2,1}-1),\qquad
\psi_{\al,n}^{(-1)}= G_{\mu_{\al,n},-1}\iv (u_{1/2,1}-1).
\ee

\begin{proposition}\label{p3.5}
Suppose (\ref{f.hal}), (\ref{f.haln}), (\ref{f.mu.an1}) and
(\ref{f.psi.pm1}).
Then (for fixed $\al>0$) the following is true:
\begin{itemize}
\item[(i)]
The operators $H(\psi_{\al,n}^{(1)})$ and $H(\psi_{\al,n}^{(-1)})$
are unitarily equivalent to the operators
$H(u_{-1/2,1})$ and $-H(u_{1/2,1})$,
respectively.
\item[(ii)]
The operators
$$
P_n(I+H(\psi_{\al,n}^{(1)}))\iv P_n-P_n
$$
are unitarily equivalent to the operators
$$
A_n=H(h_{\alpha,n})(I+H(u_{-1/2,1}))\iv H(h_{\alpha,n})-H(h_{\alpha,n})^2,
$$
which are trace class operators and converge as $n\to \iy$
in the trace norm to
$$
A=H(h_{\alpha})(I+H(u_{-1/2,1}))\iv H(h_{\alpha})-H(h_{\alpha})^2.
$$
\item[(iii)]
The operators
$$
P_n(I+H(\psi_{\al,n}^{(-1)}))\iv P_n-P_n\quad
$$
are unitarily equivalent to the operators
$$
B_n=H(h_{\alpha,n})(I-H(u_{1/2,1}))\iv H(h_{\alpha,n})-H(h_{\alpha,n})^2,
$$
which are trace class operators and converge as $n\to\iy$
in the trace norm to
$$
B=H(h_{\alpha})(I-H(u_{1/2,1}))\iv H(h_{\alpha})-H(h_{\alpha})^2.
$$
\end{itemize}
\end{proposition}
\begin{proof}
(i):\
We employ the Lemma \ref{l2.9} in order to conclude that
$$
H(\psi_{\al,n}^{(1)})= R_{\mu_{\al,n},1}^*H(u_{-1/2,1})R_{\mu_{\al,n},1},\qquad
H(\psi_{\al,n}^{(-1)})= -R_{\mu_{\al,n},-1}^*H(u_{1/2,1})R_{\mu_{\al,n},-1}.
$$

(ii):\
We first introduce the operator $W_n=H(t^n)$ and remark that $W_n^2=P_n$
and $W_nP_n=P_nW_n=W_n$. It is easily seen that the operator
$P_n(I+H(\psi_{\al,n}^{(1)}))\iv P_n-P_n$ is unitarily equivalent to
the operator $W_n(I+H(\psi_{\al,n}^{(1)}))\iv W_n-W_n^2$ by means of the
unitary and selfadjoint operator $W_n+(I-P_n)$.

Now we use the unitary equivalence established in (i) in connection with the fact that
$R_{\mu_{\al,n},1}W_nR_{\mu_{\al,n},1}^*= R_{\mu_{\al,n},1}H(t^n)R_{\mu_{\al,n},1}^*=H(h_{\al,n})$ (see again Lemma  \ref{l2.9}).
Notice that $h_{\al,n}=G_{\mu_{\al,n},1}(t^n)$.
This implies the unitary equivalence to $A_n$.

In order to prove the convergence $A_n\to A$ in the trace norm we
write
$$
A_n =H(h_{\alpha,n})(I+H(u_{-1/2,1}))\iv H(u_{-1/2,1})H(h_{\alpha,n}).
$$
The function $h_{\al,n}$ is uniformly bounded and converges (along with all its
derivatives) uniformly on each compact subset of $\T\setminus\{-1\}$ to
the function $h_\alpha$. Hence (by Lemma \ref{l2.7})
$$
H(h_{\alpha,n})\to H(h_\al),\qquad T(\widetilde{h_{\al,n}})\to T(\widetilde{h_\al})
$$
strongly on $H^2(\T)$. The same holds for their adjoints.

Next we claim that all operators
$H(u_{-1/2,1})H(h_{\alpha,n})$ are trace class operators and converge in the
trace norm to $H(u_{-1/2,1})H(h_{\alpha})$. To see this we choose two smooth
functions $f$ and $g$ on $\T$ which vanish identically in a neighborhood
of $-1$ and $1$, respectively, such that $f+g=1$. Then
we decompose
\bqn
H(u_{-1/2,1})H(h_{\alpha,n})
&=&
H(u_{-1/2,1})T(f)H(h_{\alpha,n})+
H(u_{-1/2,1})T(g)H(h_{\alpha,n})\nn\\[1ex]
&=&
H(u_{-1/2,1})H(fh_{\alpha,n})-
H(u_{-1/2,1})H(f)T(\widetilde{h_{\alpha,n}})\nn\\[1ex]
&&+
H(u_{-1/2,1}\tilde{g})H(h_{\alpha,n})-
T(u_{-1/2,1})H(\tilde{g})H(h_{\alpha,n}).\nn
\eqn
The Hankel operators $H(f)$ and $H(\tilde{g})$ are both trace class
and so are the operators $H(fh_{\alpha,n})$ and $H(u_{-1/2,1}\tilde{g})$
since the generating functions are smooth.

Moreover, $fh_{\al,n}\to fh_\al$ uniformly and the same holds for the derivatives.
Hence $H(fh_{\al,n})\to H(fh_\al)$ in the trace norm by Lemma \ref{l2.6}.
Along with the strong convergence noted above, it follows that
$H(u_{-1/2,1})H(h_{\alpha,n})$ converges in the trace norm to
$$
H(u_{-1/2,1})H(fh_{\alpha})-
H(u_{-1/2,1})H(f)T(\widetilde{h_{\alpha}})+
H(u_{-1/2,1}\tilde{g})H(h_{\alpha})-
T(u_{-1/2,1})H(\tilde{g})H(h_{\alpha}),
$$
which is trace class and equal to $H(u_{-1/2,1})H(h_{\alpha})$.

(iii):\ The proof of these assertions is analogous. The only (slight) difference is
that $R_{\mu_{\al,n},-1}W_nR_{\mu_{\al,n},-1}^*=R_{\mu_{\al,n},-1}H(t^n)R_{\mu_{\al,n},-1}^*=(-1)^{n+1}H(h_{\al,n})$
as $G_{\mu_{\al,n},-1}(t^n)=(-1)^n h_{\al,n}$.
The possibly different sign at this place does not effect the argumentation.\end{proof}


\section{Proof of the asymptotic formula}
\label{sec:5}

In this section we are going to prove the asymptotic formula
(\ref{f.Dys}).

Our first goal is to discretize the Wiener-Hopf operator $I-K_\al$,
which will lead us to a Toeplitz operator. Here and in what follows
$\chi_{\alpha}$ stands for the characteristic function of the subarc
$\{e^{i\theta}:\al<\theta<2\pi-\al\}$ of $\T$.

\begin{proposition}\label{p3.1}
For each $\al>0$ we have
\bqn
\det(I-K_\al)=\lim_{n\to\iy} \det T_n(\chi_{\frac{\al}{n}}).
\eqn
\end{proposition}
\begin{proof}
Recall that the operator $K_\alpha$ is the integral operator on
$L^2[0,\alpha]$ with the kernel $K(x-y)$, where
$$
K(x)=\frac{\sin x}{\pi x}.
$$
Introduce the $n\times n$ matrices
$$
A_n=\left(\frac{\al}{n}K\left(\frac{\al(j-k)}{n}\right)
\right)_{j,k=0}^{n-1},
\quad
B_n=\left(\frac{\al}{n}\int_{0}^1\int_0^1
K\left(\frac{\al(j-k+\xi-\eta)}{n}\right)\,d\xi d\eta\right)_{j,k=0}^{n-1}.
$$
By the mean value theorem the entries of $A_n-B_n$ can be estimated
uniformly by $O(n^{-2})$, whence it follows that the Hilbert-Schmidt
norm of $A_n-B_n$ is $O(n\iv)$.  Since the Hilbert-Schmidt norm of the
$n\times n$ identity matrix is $O(\sqrt{n})$, we obtain that the
trace norm of $A_n-B_n$ is $O(1/\sqrt{n})$.

The Fourier coefficients of $1-\chi_{\frac{\al}{n}}$ are
\bqn
[1-\chi_{\frac{\al}{n}}]_k &=&
\left\{\ba{cl}\displaystyle \frac{\al}{\pi n} & \mbox{ if }k=0\\[2ex]
\displaystyle
\frac{\sin(\frac{k\al}{n})}{\pi k} &\mbox{ if }k\neq0.\ea
\right. \nn
\eqn
Hence it follows that $T_n(\chi_{\frac{\al}{n}})=I_n-A_n$.
Introduce the isometry
$$
U_{\alpha,n}: \{x_k\}_{k=0}^{n-1}\in\C^n\mapsto\sqrt{\frac{n}{\al}}\;
\sum_{k=0}^{n-1} x_k\chi_{[\frac{\alpha k}{n},\frac{\alpha(k+1)}{n}]}
\in L^2[0,\alpha],
$$
and remark that
$$
U_{\alpha,n}^*: f\in L^2[0,\al]\mapsto
\left\{\sqrt{\frac{n}{\al}}\;
\int_{0}^\al f(x) \chi_{[\frac{\alpha k}{n},\frac{\alpha(k+1)}{n}]}\,dx
\right\}_{k=0}^{n-1}\in  \C^n.
$$
It can be verified straightforwardly, that $U_{\al,n}^*K_\al U_{\al,n}=B_n$.
Hence
\bqn
\det(I-K_\al) &=&
\det(I_n-U_{n,\al}^*K_\al U_{n,\al})= \det (I_n-B_n)
\nn\\
&\sim& \det(I_n-A_n)=
\det T_n(\chi_{\frac{\al}{n}})\nn
\eqn
as $n\to\iy$. This completes the proof.
\end{proof}

The following result has been established in \cite[Cor.~2.5]{BE1}.

\begin{proposition}\label{p3.2}
Let $b\in L^1[-1,1]$ and suppose that $b_0(x)=b_0(-x)$, where
$$
b_0(x)=b(x)\sqrt{\frac{1-x}{1+x}}.
$$
Then $\det H_n[b]=\det T_n(d)$ with $d(e^{i\theta})=b_0(\cos\frac{\theta}{2})$.
\end{proposition}

We use this result in order to reduce our Toeplitz determinant
$\det T_n(\chi_{\frac{\al}{n}})$ to a Hankel determinant.

\begin{proposition}\label{p3.3}
We have
\bqn
\det T_n(\chi_{\frac{\al}{n}}) &=& (\rho_{\al,n})^{n^2}\det H_n[b_{\alpha,n}],
\eqn
where
\be\label{f.rho}
b_{\alpha,n}(x)=\sqrt{\frac{1+\rho_{\alpha,n}x}{1-\rho_{\alpha,n}x}},\qquad
\rho_{\al,n} = \cos\left(\frac{\al}{2n}\right).
\ee
\end{proposition}
\begin{proof}
We apply Proposition \ref{p3.2} with
$d(e^{i\theta})=\chi_{\frac{\al}{n}}(e^{i\theta})$,
$b_0(x)=\chi_{[-\rho_{\alpha,n},\rho_{\alpha,n}]}(x)$, and
$$
b(x)=\sqrt{\frac{1+x}{1-x}}\chi_{[-\rho_{\alpha,n},\rho_{\alpha,n}]}(x).
$$
It follows that $\det T_n(\chi_{\frac{\alpha}{n}})=\det H_n[b]$. The
entries of $H_n[b]$ are the moments $[b]_{1+j+k}$, $0\le j,k\le n-1$.
A simple computation gives
$$
[b]_k = \frac{1}{\pi}\int_{-1}^1 b(x) (2x)^{k-1}\,dx=
\frac{(\rho_{\al,n})^k}{\pi}\int_{-1}^1
\sqrt{\frac{1+\rho_{\al,n}y}{1-\rho_{\al,n}y}}(2y)^{k-1}dy=
(\rho_{\al,n})^k [b_{\alpha,n}]_k.
$$
Now we call pull out certain diagonal matrices from the left and the right of
$H_n[b]$ to obtain the matrix $H_n[b_{\alpha,n}]$. The determinants of
the diagonal matrices give the factor $(\rho_{\alpha,n})^{n^2}$.
\end{proof}

In the following result we use the function
\bqn\label{f.psi}
\psi_{\al,n}(t) &=&
\left(\frac{1-\mu_{\al,n} t}{1-\mu_{\al,n} t\iv}\right)^{1/2}
\left(\frac{1+\mu_{\al,n}t\iv}{1+\mu_{\al,n}t}\right)^{1/2}
\chi(t),
\eqn
where $\chi(t)$ is given by (\ref{f.chi}) and where
\bqn\label{f.mu}
\mu_{\al,n}=\frac{1-\sqrt{1-\rho_{\al,n}^2}}{\rho_{\al,n}}
\eqn
with $\rho_{\al,n}$ given by (\ref{f.rho}).
Remark that $\mu_{\al,n}\in[0,1)$ satisfies condition
(\ref{f.mu.an1}).

\begin{proposition}\label{p3.4}
We have
\bqn
\lim_{n\to\iy} \det T_n(\chi_{\frac{\al}{n}}) &=&
e^{-\frac{\al^2}{8}}\lim_{n\to\iy}
\det\left[ P_n\Big(I+H(\psi_{\al,n})\Big)\iv P_n\right].
\eqn
\end{proposition}
\begin{proof}
We use Proposition \ref{p3.3}.
Since $\rho_{\al,n}=1-\frac{\al^2}{8n^2}+O(n^{-4})$ it is readily verified that
$(\rho_{\al,n})^{n^2} \to e^{-\frac{\al^2}{8}}$.
We obtain
$$
\lim_{n\to\iy} \det T_n(\chi_{\frac{\al}{n}}) =
e^{-\frac{\al^2}{8}}\lim_{n\to\iy}
\det H_n[b_{\al,n}].
$$

Now we employ Theorem \ref{t2.4} with
$$
c(e^{i\theta})=\sqrt{\frac{1+\rho_{\al,n}\cos\theta}{1-\rho_{\al,n}\cos\theta}}.
$$
Obviously, (since ($\rho_{\al,n}=2\mu_{\al,n}/(1+\mu_{\al,n}^2)$)
$$
c(t)=\sqrt{\frac{(1+\mu_{\al,n} t)(1+\mu_{\al,n} t\iv)}{
(1-\mu_{\al,n} t)(1-\mu_{\al,n} t\iv)}},
$$
whence we conclude that $c(t)=\tc_+(t)Gc_+(t)$ with $G=1$ and
$$
c_+(t)=
\left(\frac{1+\mu_{\al,n} t}{1-\mu_{\al,n} t}\right)^{1/2}.
$$
Furthermore,
$$
\tc_+(t)c_+\iv(t)=
\left(\frac{1-\mu_{\al,n} t}{1-\mu_{\al,n}t\iv}\right)^{1/2}
\left(\frac{1+\mu_{\al,n} t\iv}{1+\mu_{\al,n}t}\right)^{1/2}.
$$
It follows that
$$
\det H_n[b_{\al,n}] =\det\left[P_n\Big(I+H(\psi_{\al,n})\Big)\iv P_n\right].
$$
This implies the desired assertion.
\end{proof}

In the following proposition we identify the limit
of the determinant
$$
\det\left[P_n\Big(I+H(\psi_{\al,n})\Big)\iv P_n\right]
$$
as $n\to\iy$. Recall the definitions (\ref{f.u}), (\ref{f.hal}), and
Proposition \ref{p2.2}.

\begin{proposition}\label{p3.5x}
We have
\bqn
\lefteqn{
\lim_{n\to\iy}
\det \Big[P_n(I+H(\psi_{\al,n}))\iv P_n\Big]} \nn\\[1ex]
&=&
\det \Big[H(h_\al)(I+H(u_{-1/2,1}))\iv H(h_\al)\Big]
\det \Big[H(h_\al)(I-H(u_{1/2,1}))\iv H(h_\al)\Big],
\label{f.36}
\eqn
where all expressions on the right hand side are well defined.
\end{proposition}
\begin{proof}
First of all we remark that the right hand side is well defined. The
inverses exist due to Proposition \ref{p2.2}.
Notice that $H(h_\al)^2$ is a projection operator since (by (\ref{f.Tab}) and (\ref{f.THspez}))
$$
H(h_\al)^3=(I-T(h_\al)T(\tilde{h}_\al))H(h_\al)=H(h_\al).
$$
We consider the operators $H(h_\al)(I\pm H(u_{\mp 1/2,1}))\iv H(h_\al)$
as being restricted onto the image of $H(h_\al)^2$. We can complement these operators
with the projection $I-H(h_\al)^2$ without changing the value of the
corresponding determinant,
$$
\det \Big[H(h_\al)(I\pm H(u_{\mp 1/2,1}))\iv H(h_\al)\Big] =
\det \Big[I+H(h_\al)(I\pm H(u_{\mp 1/2,1}))\iv H(h_\al)-H(h_\al)^2\Big].
$$
By Proposition \ref{p3.5}(ii)-(iii) we see that this last operator
determinant is well-defined.

With $\mu=\mu_{\al,n}$ given by (\ref{f.mu}) we obtain from
(\ref{f.chi}), (\ref{f.psi}) and (\ref{f.Giv}) that
$$
\psi_{\al,n}(t) =
\left(\frac{t-\mu}{1-\mu t}\right)^{-1/2}
\left(\frac{t+\mu}{1+\mu t}\right)^{1/2}
= G_{\mu,1}\iv(u_{-1/2,1}) G_{\mu,-1}\iv(u_{1/2,1}).
$$
Introduce the functions $\psi_{\al,n}^{(\pm1)}$ by (\ref{f.psi.pm1}).
Then
$$
\psi_{\al,n} =
(\psi_{\al,n}^{(1)}+1)(\psi_{\al,n}^{(-1)}+1).
$$
Proposition  \ref{p2.8} implies that
$$
H(\psi_{\al,n}) = H(\psi_{\al,n}^{(1)})+
H(\psi_{\al,n}^{(-1)})+H(\psi_{\al,n}^{(-1)})H(\psi_{\al,n}^{(1)})
+o_1(1),
$$
where $o_1(1)$ stands for a sequence of operators converging in the trace
norm to zero as $n\to\iy$. By Proposition \ref{p3.5}(i) and Proposition \ref{p2.2},
the operators
$I+H(\psi_{\al,n}^{(1)})$ and $I+H(\psi_{\al,n}^{(-1)})$ are invertible
and their inverses are uniformly bounded. Hence
$$
(I+H(\psi_{\al,n}))\iv =
(I+H(\psi_{\al,n}^{(1)}))\iv (I+H(\psi_{\al,n}^{(-1)}))\iv +o_1(1).
$$
Using the formula $(I+A)\iv=I-(I+A)\iv A=I-A(I+A)\iv$, we can write this as
\bqn
(I+H(\psi_{\al,n}))\iv &=& -I+(I+H(\psi_{\al,n}^{(1)}))\iv+
(I+H(\psi_{\al,n}^{(-1)}))\iv\nn\\[.5ex]
&&
+(I+H(\psi_{\al,n}^{(1)}))\iv
H(\psi_{\al,n}^{(1)})H(\psi_{\al,n}^{(-1)})
(I+H(\psi_{\al,n}^{(-1)}))\iv+o_1(1).\nn
\eqn
It follows that
\bqn
P_n(I+H(\psi_{\al,n}))\iv P_n&=& -P_n+P_n(I+H(\psi_{\al,n}^{(1)}))\iv P_n+
P_n(I+H(\psi_{\al,n}^{(-1)}))\iv P_n\nn\\[.5ex]
&&
+P_n(I+H(\psi_{\al,n}^{(1)}))\iv
H(\psi_{\al,n}^{(1)})H(\psi_{\al,n}^{(-1)})
(I+H(\psi_{\al,n}^{(-1)}))\iv P_n\nn\\[.5ex]
&&+o_1(1).\nn
\eqn
Since $I-P_n=I-H(t^n)^2=T(t^n)T(t^{-n})$ (see (\ref{f.Tab})), we have
\bqn
H(\psi_{\al,n}^{(1)})(I-P_n)H(\psi_{\al,n}^{(-1)}) &=&
H(\psi_{\al,n}^{(1)})T(t^n)T(t^{-n})H(\psi_{\al,n}^{(-1)})\nn\\
&=&
T(t^{-n})H(\psi_{\al,n}^{(1)})H(\psi_{\al,n}^{(-1)})T(t^n)
=o_1(1),\nn
\eqn
where we used also Proposition \ref{p2.8} and (\ref{f.THspez}).
Hence we obtain
\bqn
P_n(I+H(\psi_{\al,n}))\iv P_n&=& -P_n+P_n(I+H(\psi_{\al,n}^{(1)}))\iv P_n+
P_n(I+H(\psi_{\al,n}^{(-1)}))\iv P_n\nn\\[.5ex]
&&
+P_n(I+H(\psi_{\al,n}^{(1)}))\iv
H(\psi_{\al,n}^{(1)}) P_n H(\psi_{\al,n}^{(-1)})
(I+H(\psi_{\al,n}^{(-1)}))\iv P_n\nn\\[.5ex]
&&+o_1(1)\nn\\[1ex]
&=&
P_n(I+H(\psi_{\al,n}^{(1)}))\iv P_n (I+H(\psi_{\al,n}^{(-1)}))\iv P_n+o_1(1).\nn
\eqn
Therein we employed again the formula $(I+A)\iv=I-(I+A)\iv A=I-A(I+A)\iv$.
Proposition \ref{p3.5}(ii)-(iii) implies the
uniform boundedness of the operators $P_n(I+H(\psi_{\al,n}^{(1)}))\iv P_n$
and $ P_n (I+H(\psi_{\al,n}^{(-1)}))\iv P_n$. In connection with the
well-known formula
$$
|\det(I+A)-\det(I+B)|\;\le\;\|A-B\|_1\exp(\max\{\|A\|_1,\|B\|_1\}),
$$
this proves that
\bqn
\lefteqn{\lim_{n\to\iy}
\det\Big[P_n(I+H(\psi_{\al,n}))\iv P_n\Big]} \nn\\
&=&
\lim_{n\to\iy}
\det \Big[P_n(I+H(\psi_{\al,n}^{(1)}))\iv P_n\Big]
\det\Big[P_n(I+H(\psi_{\al,n}^{(-1)}))\iv P_n\Big].\nn
\eqn
These determinants can be written as
$$
\det \Big[P_n(I+H(\psi_{\al,n}^{(\pm1)}))\iv P_n\Big]=
\det \Big[I+P_n(I+H(\psi_{\al,n}^{(\pm1)}))\iv P_n-P_n\Big],
$$
and now the convergence in the trace norm stated in Proposition
\ref{p3.5}(ii)-(iii) implies the desired assertion.
We remark in this connection that
\bqn\label{f.mu-exp}
\mu_{\al,n}=1-\frac{\al}{2n}+O(n^{-2}),\qquad n\to\iy,
\eqn
holds.
\end{proof}

In regard to the next result, recall the definition (\ref{f.uhat}) of the
functions $\hat{u}_\beta$.

\begin{theorem}
We have
\bqn
\det(I-K_\al) &=& \exp\left(-\frac{\al^2}{8}\right)
\det\Big[\Pi_{\frac{\al}{2}}(I+H_{\R}(\hat{u}_{-1/2}))\iv \Pi_{\frac{\al}{2}})\Big]
\nn\\[1ex]
&&\times
\det\Big[\Pi_{\frac{\al}{2}}(I-H_{\R}(\hat{u}_{1/2}))\iv \Pi_{\frac{\al}{2}}\Big],
\label{f.64}
\eqn
where all expressions on the right hand side are well defined.
\end{theorem}
\begin{proof}
We combine Proposition \ref{p3.1} with Proposition \ref{p3.4}
and Proposition \ref{p3.5x} to conclude that
\bqn
\det(I-K_\al) &=& \exp\left(-\frac{\al^2}{8}\right)
\det \Big[H(h_\al)(I+H(u_{-1/2,1}))\iv H(h_\al)\Big]\nn\\
&&\times
\det \Big[H(h_\al)(I-H(u_{1/2,1}))\iv H(h_\al)\Big].\nn
\eqn
By means of the transform $S$ considered in (\ref{f.cd}) and (\ref{f.frac})
we notice that
$$
H_\R(\hat{u}_\beta)=SH(u_{\beta,1})S^*,\qquad
H_\R( e^{ix\alpha/2})=S H(h_\al) S^*.
$$
It remains to remark that $H( e^{ix\alpha/2})^2=\Pi_{\al/2}$.
\end{proof}

Now we use the following asymptotic formula for the two operator determinants
appearing on the right hand side of (\ref{f.64}), which is proved in \cite{BE3}.
Therein $G(z)$ stands for the Barnes $G$-function \cite{Bar}.
For convenience we make a change in variables $\al\mapsto 2\al$.

\begin{theorem}
The following asymptotic formulas hold,
\bqn
\det\Big[\Pi_{\al}(I+H_{\R}(\hat{u}_{-1/2}))\iv \Pi_{\al}\Big]
&\sim&
\al^{-1/8}\pi^{1/4} 2^{1/4}G(1/2),\qquad \al\to\iy,\\
\det\Big[\Pi_{\al}(I-H_{\R}(\hat{u}_{1/2}))\iv \Pi_{\al}\Big]
&\sim&
\al^{-1/8}\pi^{1/4} 2^{-1/4}G(1/2),\qquad \al\to\iy.
\eqn
\end{theorem}

Combining the previous results we get the desired asymptotic formula.

\begin{theorem}
The asymptotic formula
\bqn\label{f.as1}
\log\det(I-K_{2\al}) &=& -\frac{\al^2}{2}
-\frac{\log\al}{4}+C+o(1),\qquad
\al\to\iy,
\eqn
holds with the constant
\bqn\label{f.con2}
C=\frac{\log 2}{12}+3\zeta'(-1).
\eqn
\end{theorem}
\begin{proof}
The previous two theorems give the asymptotic formula
\bqn
\det(I-K_{2\al}) &\sim& \exp\left(-\frac{\al^2}{2}\right)
\al^{-1/4}\pi^{1/2}(G(1/2))^2,\qquad \al\to\iy.
\eqn
We can express $G(1/2)$ in terms of $\zeta'(-1)$, where $\zeta$
is Riemann's zeta function.
According to \cite[page 290]{Bar} we have
\bqn
\log G(1/2) = -\frac{\log\pi}{4} +\frac{1}{8}-\frac{3\log A}{2}+\frac{\log 2}{24}
\nn
\eqn
with $A= \exp(-\zeta'(-1)+1/12)$ being Glaisher's constant.
Hence
\bqn
2\log G(1/2) = -\frac{\log\pi}{2}+3\zeta'(-1)+\frac{\log 2}{12},\nn
\eqn
which implies the desired asymptotic formula (\ref{f.as1}) with the
constant (\ref{f.con2}).
\end{proof}



\end{document}